\documentclass[11pt]{article}

\usepackage[a4paper,margin=2.5cm]{geometry}

\usepackage{newtxtext,newtxmath}

\usepackage[english]{babel}

\usepackage{amsmath,bm}
\usepackage{setspace}
\usepackage{graphicx}
\usepackage{booktabs}
\usepackage{float}
\usepackage{caption}
\usepackage[section]{placeins}

\usepackage{tikz}
\usetikzlibrary{arrows.meta}
\tikzset{>={Latex[width=2mm,length=2mm]}}

\usepackage{authblk}

\usepackage{xurl}
\usepackage[hidelinks]{hyperref}

\usepackage[
  backend= bibtex, %biber,
  style=authoryear,
  sorting=nyt,
  maxbibnames=99,
  giveninits=true,
  doi=true,
  url=false,
  isbn=false
]{biblatex}

\DeclareFieldFormat{doi}{%
  \href{https://doi.org/#1}{https://doi.org/#1}}

\newcommand\Bigskip{\bigskip\noindent}

\newcommand\BEQ{\begin{equation}}
\newcommand\EEQ{\end{equation}}
\newcommand\BEA{\begin{eqnarray}}
\newcommand\EEA{\end{eqnarray}}

\begin{document}

\title{\Large\bfseries Network analysis reveals phase transitions in agro-food nitrogen systems with contrasting feeding capacity and land requirement}

\author[1,2]{Adrien Fauste-Gay\thanks{Corresponding author: adrien.fauste-gay@univ-grenoble-alpes.fr}}
\author[3]{Jérémie Unterberger}
\author[1]{Olivier Vidal}
\author[2]{Lauriane Mouysset}

\affil[1]{CNRS, ISTerre, University of Grenoble Alpes, Grenoble, France}
\affil[2]{CNRS, CIRED, Nogent-sur-Marne, France}
\affil[3]{IECL, Université de Lorraine, Vandœuvre-lès-Nancy, France}

\date{}

\maketitle
%\linenumbers

\vspace{2mm}
\begin{abstract}

Agro-food transitions are commonly assessed using optimization models or scenario approaches that return a single feasible configuration. These approaches provide estimates of system performance but limited information on feasible nitrogen cycling configurations. In this study, we develop a network-based modelling framework to characterize alternative steady-state configurations of agro-food nitrogen systems and identify phase transitions between these configurations.

We formulate a continuous-time compartmental network model in which production, allocation, recycling, and losses are represented as mass-conserving flows. The system is analyzed under steady-state conditions. A coarse-graining procedure identifies dominant recycling cycles and derives aggregate indicators describing system configuration. The model is parameterized using French reference data and applied to case studies covering fertilization, livestock, and dietary patterns.

The analysis identifies distinct steady-state configurations of the nitrogen network, separated by phase transitions in dominant recycling structures. These transitions correspond to changes in system performance indicators. For the French reference system, cropland-based feeding capacity is about 8.7 people per hectare of cropland (depending on constraints), while strict dietary autonomy requires about 19.7 million hectares of agricultural land. Reducing synthetic fertilization without changes in diet or nitrogen sourcing leads to higher land requirements or greater dependence on external inputs.

The framework provides a method for mapping feasible nitrogen system configurations and analyzing phase transitions in agro-food nitrogen networks under alternative parameters and constraints.

\end{abstract}

\noindent\textbf{Keywords:} agro-food system transition; nitrogen metabolism; network analysis; feeding capacity; land requirement; crop-livestock integration; dietary transition; agricultural autonomy; synthetic fertilization; biological nitrogen fixation

%\keywords{agro-food system transition; nitrogen metabolism; network analysis; feeding capacity; land requirement; crop-livestock integration; dietary transition; agricultural autonomy; synthetic fertilization; biological nitrogen fixation}

%\tableofcontents

%\newpage

\section{Introduction}

Agricultural transitions are largely constrained by the way reactive nitrogen is produced, allocated, recycled, and lost across cropland, grassland, livestock, human consumption, and the environment. Because nitrogen flows couple yields, feed requirements, dietary supply, and pollution, they are often used in agricultural research to characterise agricultural systems. Moreover, forecast analysis must satisfy mass balance and agronomic bounds while remaining explicit enough to identify actionable levers (e.g., manure routing, grazing share, fertilizer substitution, and trade dependence). This requirement is particularly acute in a context of climate and land-use pressures \parencite{IPCC2019Land}, biodiversity loss \parencite{Tilman2002}, nutrient pollution \parencite{conley2009controlling}, and market volatility \parencite{dupas2019time}, which jointly constrain feasible configurations of agro-food systems.

Socio-metabolic and metabolic flux analysis (MFA) approaches provide an appropriate biophysical lens by representing agro-food systems as networks of material flows, where conservation laws, technical coefficients, and admissible couplings are explicit \parencite{fischerKowalski2015social, haberl2019contributions}. In the case of agriculture, this ``metabolic'' perspective makes it possible to track nitrogen circulation between production compartments and final uses, and to distinguish what is physically feasible from what is merely selected by a particular objective or narrative. This paper builds on that representation to move from single-solution allocations toward a systematic exploration of feasible operating modes.

GRAFS is a model built on this approach: it represents territorial agro-food metabolism with 4 agricultural compartments (Population, Cropland, Grassland, Livestock) and boundary compartments (trade, Haber-Bosch, Environment). In practice, GRAFS allows to analyze functioning of past and present agricultural systems in various contexts while keeping the accounting physically coherent \parencite{le2018biogeochemical}.

For forecast, optimization-based MFA tools provide one feasible allocation under a chosen objective, but do not characterize the diversity of allocations compatible with the same physical constraints. Here we recast agro-food nitrogen metabolism in a formalism which allows analytical exploration of the feasible space and its qualitatively distinct operating modes, called {\em metabolic phases}. We show how these phases translate into contrasted outcomes for sustainable population density, self-sufficiency, and dependence on industrial nitrogen inputs.

More broadly, agricultural foresight still oscillates between coarse aggregate indicators and general narratives. What is lacking are tools that systematically explore the full space of physically possible configurations, the set of allocations consistent with conservation, agronomic bounds, and demand—so stakeholders can debate realistic futures on the basis of what could work, not just what an objective function selects. Optimization remains useful to select a particular configuration; our goal is to characterize the set of configurations consistent with the same physical constraints.

In this paper, we propose a general mathematical framework to map and explore the feasible metabolic solution space, thereby linking physically realistic forecast analysis to actionable levers for transforming agro-food systems. The starting point is a reformulation of the model in terms of {\em stock variables} instead of {\em flux variables}. {\em Stock variables} represent the part of the nitrogen stock of compartments which is available for biophysical or economical exchange. They are commonly used in more complex agent-based models like MAELIA that discuss agricultural land parcel dynamics, and allow in particular such features as  a more general discussion of slow chemical transformations of organic nitrogen in the soil delaying its availability for plants, and  non-stationary situations \parencite{gaudou2014maelia}. 

We use here an intermediate formalism originated from the study of {\em open chemical reaction networks}. There is no spatialization so that compartment nitrogen stocks may be thought as interacting chemical species in solution transforming one into the other through chemical mechanisms. Though the reformulation is straightforward for biophysical processes, this should be thought on the whole as a formal analogy, allowing a transfer of mathematical techniques used in the field of kinetic reaction networks. We build on recent methodological results \parencite{UntNgh,Unt-math,Unt-phys} allowing a coarse-grained description of the dynamics of complex chemical mixtures, that are surprisingly well-suited for a systematic exploration of the  {\em phases} of nitrogen metabolism.

{\em Outline of the paper.} We start with an accurate description of the network used in this paper (fluxes, stocks) and of its main variables of interest and provide a short summary of the mathematical methodology in the Materials and methods section (Section \ref{section:methods}); the interested reader may consult the Supplementary Materials (Suppl. Mat.) file for detailed explanations. The Results section (Section \ref{section:results}) describes our findings.

%%%%%%%%%%%%%%%%%%%%%%%%
%%%%%%%%%%%%%%%%%%%%%%ù

\section{Materials and methods}
\label{section:methods}

%%%%%%%%%%%%%%%%%%%%
%%%%%%%%%%%%%%%%%%%%%%%%%%%%

We analyze agro-food system dynamics using a network-based representation derived from the GRAFS formalism. To illustrate the proposed methodology, we introduce a deliberately simplified “toy” nitrogen metabolism network that captures the core structural features of territorial agro-food systems while remaining analytically tractable.

\begin{figure}[tbp]
    \centering
    \includegraphics[width=0.8\textwidth]{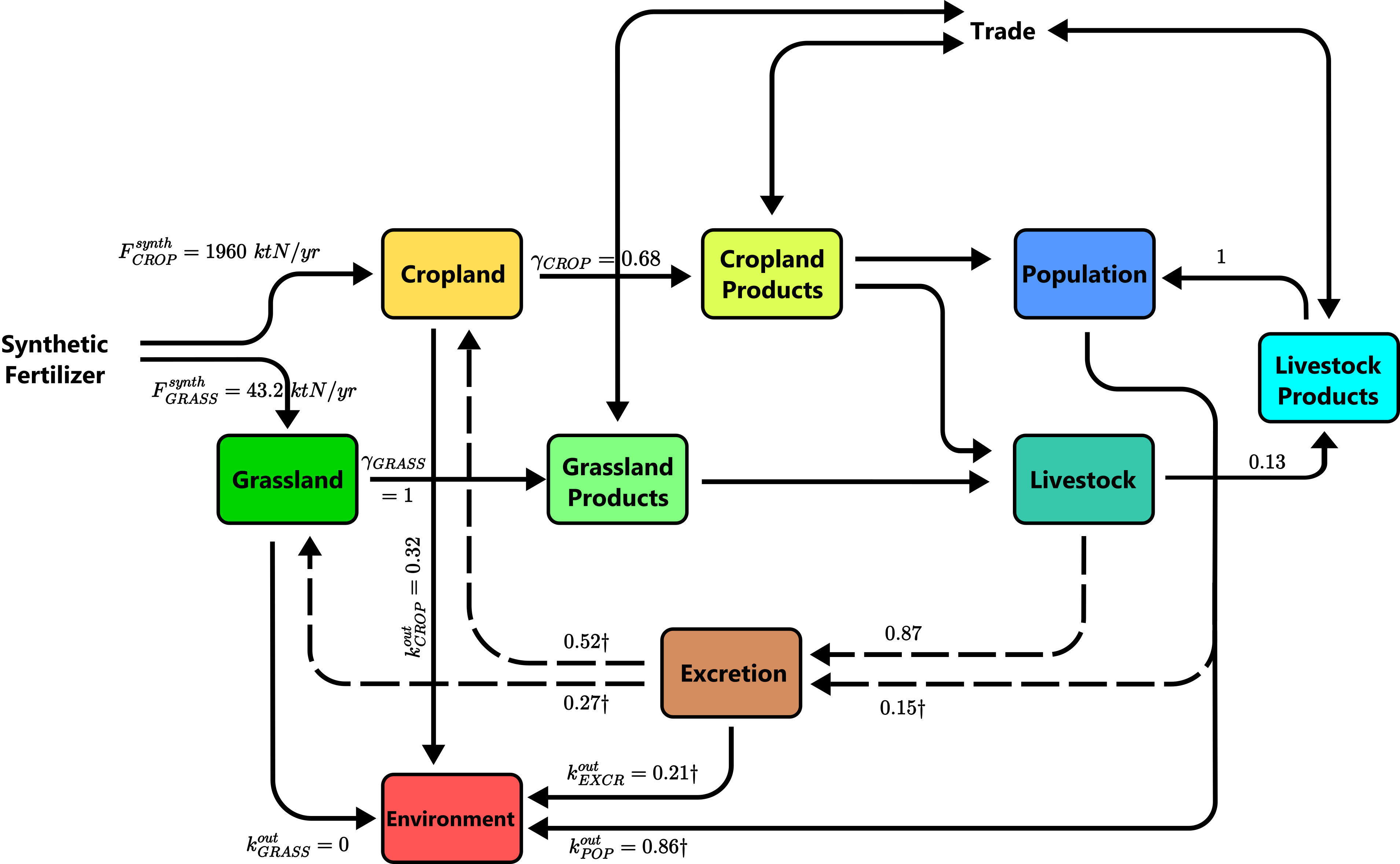}
    \caption{Stylized nitrogen network used in this study in the unsaturated regime. Nodes represent cropland, grassland, agricultural products, livestock, population, excretion, and environmental losses; dashed arrows indicate recycling flows, and the values shown are 2014 French reference coefficients.}
    \label{fig:G0}
\end{figure}

The network consists of nine compartments: cropland, grassland; three agricultural product compartments (cropland products, grassland products, and livestock products); livestock, population; excretion; and, finally, the environment (atmosphere, groundwater and soil). This structure is directly adapted from the original GRAFS model \parencite{le2018biogeochemical}, which preceded the optimization-based E-GRAFS framework. The resulting network topology is shown in Figure~\ref{fig:G0}.

We use abbreviated notations for compartments: CROP (cropland) and GRASS (grassland); CRP (cropland products) and GRP (grassland products);
POP (human population) and LVS (livestock); LVP (livestock products); EXC (excretion); and $\emptyset$ (environment).

\bigskip {\bf Fluxes.} 
Network compartments are connected by 14 nitrogen fluxes, denoted by $F$ and  expressed in ktN/yr. Production fluxes connect land compartments to product compartments: $F_{CROP \rightarrow CRP}$ (crop production) and $F_{GRASS \rightarrow GRP}$ (grassland production). Cropland products are allocated either directly to the human population ($F_{CRP \rightarrow POP}$, vegetal food) or to livestock ($F_{CRP \rightarrow LVS}$), while grassland products are exclusively used as livestock feed ($F_{GRP \rightarrow LVS}$). Livestock production is represented by the flux $F_{LVS \rightarrow LVP}$, and livestock products are consumed by the population through $F_{LVP \rightarrow POP}$.

\medskip
Nitrogen recycling is represented by return fluxes (shown as dashed arrows in the network): excretion from humans and livestock ($F_{POP \rightarrow EXC}$ and $F_{LVS \rightarrow EXC}$), followed by the redistribution of excreta to agricultural soils ($F_{EXC \rightarrow CROP}$ and $F_{EXC \rightarrow GRASS}$).

\medskip
Nitrogen losses to the environment are modeled as volatilization or leakage fluxes from cropland, grassland, and excretion compartments ($F_{CROP \rightarrow \emptyset}$, $F_{GRASS \rightarrow \emptyset}$, $F_{EXC \rightarrow \emptyset}$). Since there are no return flows from the environment, it is natural to  consider it  to be an external compartment formally denoted 
$\emptyset$; then losses may be seen as outflows from the system, and denoted for
short $F^{out}_{CROP}, F^{out}_{GRASS}, F^{out}_{EXC}$.

\medskip These {\em internal fluxes} are complemented by {\em external fluxes}: synthetic nitrogen inputs $F^{synth}_{CROP}, F^{synth}_{GRASS}$ and biological nitrogen fixation fluxes $BNF_{CROP},BNF_{GRASS}$, as opposed to net imports or exports of tradable agricultural products,
\begin{equation}
\boldsymbol{\Phi}^{prod} = (\Phi_{CRP}, \Phi_{GRP}, \Phi_{LVP}),
\end{equation}

%\bigskip {\bf Stocks.} {\color{red} In this simplified model, all nitrogen taken in by an internal compartment in the course of a year is either lost to the environment or redistributed to other internal compartments through the fluxes $f$. The approximation overlooks chemical transformations and transport phenomena taking place over a period of time spanning several years, like (de)nitrification, and also assumes stationarity. As a result, each compartment $i$ contains an equivalent circulating stock $v_i$ equal to the sum of fluxes directed to (equivalently, from) that compartment,  $v_i = \sum_j F_{j\to i} = \sum_j F_{i\to j}$.  All  internal fluxes -- those of biophysical origin, but also allocation fluxes -- are in a natural way a fraction of stocks, whereas external fluxes are mainly dependent on agriculture management.}

\bigskip {\bf Stocks.} In this simplified representation, nitrogen entering a compartment over a year follows two main fates: it is either lost to the environment (e.g., volatilization and other diffuse losses) or transferred to other compartments through production, feeding, consumption, and recycling flows. The model is deliberately annual and stationary: it does not resolve multi-year carry-over and slow processes such as soil organic matter turnover, delayed mineralization, or long-range transport and storage that can decouple inputs from outputs over several years. 

Under these assumptions, each compartment is characterized by an equivalent \emph{circulating stock} 
\BEQ v_i = \sum_j F_{i\to j}, \label{eq:circulating-stock}
\EEQ
 interpreted as the annual nitrogen throughput of that compartment, i.e.\ the total amount of nitrogen transiting through it in a year. Internal transfers between compartments are naturally expressed as fractions of these throughputs ({\em allocation rates}), whereas external inputs and outputs (synthetic fertilizer, biological fixation, trade, and losses) are primarily shaped by management choices and infrastructure.

%%%%%%%%%%%%%%%%%%%%%%%%%%%%%%%%

\subsection{Main variables of interest}
\label{subsection:main-variables}

%%%%%%%%%%%%%%%%%%%%%%%%%%ù

Before we describe in detail the technical parameters of the model, we emphasize the five main
variables. All numerical data are extracted or inferred from Le Noë's thesis and associated data \parencite{le2018biogeochemical}, see Suppl. Mat. for more details.

\medskip

\begin{itemize}
\item[(1)] \textsc{VolRate} $= k^{out}_{EXC}$, a technical parameter associated with nitrogen losses during
excretion and manure management.
In the current French context, \textsc{VolRate} is approximately $21\%$, but it may decrease
under extensive grazing systems (down to about $10\%$ for full grazing) or increase under
confined livestock and barn-based manure management;

\item[(2)] \textsc{RecRate}, which reflects the efficiency of nitrogen recovery from human
excreta and is strongly influenced by wastewater treatment infrastructure and lifestyles.
Historically, \textsc{RecRate} reached values close to $90\%$ in traditional rural China,
whereas it remains below $15\%$ in contemporary France (urban and rural combined)
\parencite{starck2023excreta};

\item[(3)]  the total utilized agricultural area ${\cal S}^{tot}$, and its allocation
between cropping (${\cal S}_{CROP}$) and pasture (${\cal S}_{GRASS}$);

\item[(4)] the allocation parameter $k_{CRP \to LVS}$, which governs the share of cropland
products directed toward livestock feeding.
Within an input-driven metabolic network, this parameter is formally free; in practice,
imposing zero net imports of cropland products allows it to be fixed \emph{a posteriori}
for a given livestock population and prescribed diet;

\item[(5)] finally, the allocation of excreta between cropland and grassland, characterized by the
ratio $k_{EXC \to CROP}/k_{EXC \to GRASS}$.
This ratio depends on livestock management practices as well as the relative availability
of cropland and grassland.
In present-day France, excreta allocation per ha is approximately balanced, with
\BEQ
\frac{k_{EXC \to CROP}}{k_{EXC \to GRASS}}
= \frac{0.66}{0.34}
\simeq 2 \simeq \frac{{\cal S}_{CROP}}{{\cal S}_{GRASS}}   \label{eq:EAB}
\EEQ

\end{itemize}

\medskip
\noindent\emph{Controllability of parameters.}
The parameters listed below differ markedly in their degree of \emph{pilotability}, i.e.\ the extent to which they can be modified by deliberate human action versus being constrained by biophysical conditions.
Some quantities are primarily biophysical (e.g.\ yield-response parameters such as $\alpha$ and saturation thresholds), yet they may evolve over longer timescales through plant breeding, changes in crop rotations, and shifts in pedoclimatic conditions.
Other parameters are strongly constrained by biophysics but remain adjustable within a practical range through management and infrastructure choices - notably excretion losses and collection (\textsc{VolRate}, \textsc{CollRate}), the share of grazing versus housing, and allocation decisions along the food chain (e.g.\ $k_{CRP\to LVS}$ and the partition of excreta between cropland and grassland).
Finally, a subset of variables is essentially free within the network formalism and is therefore best interpreted as policy- or behavior-driven degrees of freedom, whose values must be fixed by additional constraints or scenario assumptions.

\bigskip As we shall see, these variables are the main ones driving phase differentiation.
 We now discuss technical coefficients in detail. 

%%%%%%%%%%%%%%%%%%%%%%%%%%%%ù
%%%%%%%%%%%%%%%%%%%%%%%

\subsection {External fluxes}

The CROP and GRASS compartments represent nitrogen stocks in cropland and grassland soils. They are extensive quantities, proportional to the corresponding surface areas
${\cal S}_{CROP}$ and ${\cal S}_{GRASS}$.

In addition to stock-dependent internal fluxes, we introduce external nitrogen inputs to agricultural soils in the form of synthetic fertilizer applications and biological nitrogen fixation (BNF). Synthetic inputs are represented by constant fluxes $F^{synth}_{CROP}$ and $F^{synth}_{GRASS}$. Biological nitrogen fixation is modeled as
\begin{equation}
BNF_i = C_{T,i} \times {\cal S}_i \times \tau_{LEG,i},
\qquad i \in \{CROP, GRASS\},
\end{equation}
where $\tau_{LEG,i}$ denotes the share of legumes (by convention equal to~1 for
$ i = GRASS$), and $C_{T,i}$ (in kgN\,ha$^{-1}$) is a technical coefficient.

For convenience, we introduce intensive (per-hectare) fluxes, denoted by lowercase symbols:
\begin{equation}
\begin{aligned}
f^{synth}_{CROP} &= \frac{F^{synth}_{CROP}}{{\cal S}_{CROP}}, \qquad
f^{synth}_{GRASS} = \frac{F^{synth}_{GRASS}}{{\cal S}_{GRASS}}, \\ \\
\mathrm{bnf}_{CROP} &= C_{T,CROP}\times \tau_{LEG,CROP}, \qquad
\mathrm{bnf}_{GRASS} = C_{T,GRASS}.
\end{aligned}
\label{eq:intensive-fluxes}
\end{equation}

The external nitrogen input to each land compartment is defined as the sum of synthetic and biological inputs:
\begin{equation}
F_i = F^{synth}_i + BNF_i,
\qquad i \in \{CROP, GRASS\}.
\end{equation}
Dividing by surface area yields the corresponding intensive inputs:
\begin{equation}
f_i = f^{synth}_i + \mathrm{bnf}_i,
\qquad i \in \{CROP, GRASS\}.
\end{equation}

%%%%%%%%%%%%%%%%%%%%%%%%%
 
 \subsection{Collecting rate parameter and human excretion losses}
 \label{subsection:CollRate}

Human and livestock excreta differ in their effective loss fractions to the environment:
human excreta typically undergo much larger losses (about $90\%$) \parencite{starck2023excreta},
whereas livestock excreta losses are substantially lower (about $32\%$).
Representing these two loss processes with two distinct excretion compartments would be
straightforward but would increase network complexity. We therefore keep a single
excretion compartment \textsc{EXC} and account for the difference through the
population outflows.

We decompose the outflow from \textsc{POP} into two branches:
(i) a \emph{collected} fraction that enters \textsc{EXC}, and (ii) a \emph{direct loss}
fraction that goes immediately to the environment:
\[
k_{POP\to EXC} + k_{POP\to \emptyset} = 1,
\qquad
\textsc{CollRate} := k_{POP\to EXC}.
\]
Thus \textsc{CollRate} is the collection (or capture) rate of human excreta into the common
excretion compartment.

Within \textsc{EXC}, we use a single loss parameter
\[
\textsc{VolRate} := k_{EXC\to \emptyset},
\]
chosen to represent livestock-related excretion losses. The remaining fraction
$(1-\textsc{VolRate})$ is available for redistribution to agricultural soils.

Under this convention, the \emph{effective} recycling fraction of human excretion
(i.e.\ the share that is collected \emph{and} not lost from \textsc{EXC}) is
\begin{equation}
\textsc{RecRate}
:= k_{POP\to EXC}\,(1-k_{EXC\to \emptyset})
= \textsc{CollRate}\,(1-\textsc{VolRate}).
\label{eq:recrate-definition}
\end{equation}
Equivalently, for a prescribed recycling target \textsc{RecRate} and a given \textsc{VolRate},
the required collection rate is
\begin{equation}
\textsc{CollRate}
= \frac{\textsc{RecRate}}{1-\textsc{VolRate}}.
\label{eq:collrate-recrate}
\end{equation}
This construction captures the higher loss rate of human excreta while retaining a single
excretion compartment and a unique excretion-to-environment allocation parameter.

%%%%%%%%%%%%%%%%%%%%%%%%%%%ù

\subsection{Internal fluxes}

%%%%%%%%%%%%%%%%%%%%%%%%%%%

\subsubsection{Direct and reverse allocation ratios}

In a broader sense, the {\bf allocation ratio} of the N stock from compartment $i$ into compartment $j$ is 
the ratio $\frac{F_{i\to j}}{\sum_j F_{i\to j}} = \frac{F_{i\to j}}{v_i} \in [0,1]$, see (\ref{eq:circulating-stock}).  Some of these ratios are mainly of
biophysical origin (yields, volatilization rate \textsc{VolRate}; this also includes the recycling rate \textsc{RecRate}), and are discussed separately in \S \ref{subsection:CollRate}--\ref{subsection:yields}. We discuss here allocation ratios in a narrower sense, those
which are mainly influenced by agriculture management. There are only three in this model, which are directly
associated with variables (4) and (5) in \S \ref{subsection:main-variables}  : 

\begin{itemize}
    \item[(i)] the  {\em crop allocation ratio}, defined either as $k_{CRP\to POP}$ or  $k_{CRP\to LVS} = 1-k_{CRP\to POP}$.  It is not fixed in the model. However, it is constrained by the livestock diet
    (see (\ref{eq:diet}) below). 
    \item[(ii)] the {\em excreta allocation ratios}, defined  as $k_{EXC\to CROP}, k_{EXC\to GRASS}$, which
    do not sum up to 1 but to $1-k^{out}_{EXC}$.
\end{itemize}

Allocation ratios (except (i)) are reported on  Fig. \ref{fig:G0}; they are emphasized by a dagger ($\dagger$).
 
 \bigskip These allocation ratios are called {\em direct}, as opposed to {\bf reverse allocation ratios}, 
 which describe how a target compartment draws nitrogen from multiple
sources. Ratios characterizing {\bf human and livestock diets},
\BEQ  (\alpha_{POP},1-\alpha_{POP})=
\left(
\frac{f_{CRP \to POP}}{v_{POP}},
\frac{f_{LVP \to POP}}{v_{POP}}
\right),
\qquad (\alpha_{LVS},1-\alpha_{LVS})=
\left(
\frac{f_{CRP \to LVS}}{v_{LVS}},
\frac{f_{GRP \to LVS}}{v_{LVS}}
\right)
\label{eq:diet} 
\EEQ
are important decision variables.

%%%%%%%%%%%%%%%%%%%%
\subsubsection{Yields}  \label{subsection:yields}
%%%%%%%%%%%%%%%%%%%%

Based on values reported in the literature (see Suppl. Mat.), crop and grassland yields are modeled as functions of nitrogen inputs using a linear response that saturates beyond a given threshold \parencite{cerrato1990comparison}. This formulation captures a key agronomic characteristic of nitrogen application, namely a strictly positive response up to a saturation threshold, followed by a zero marginal response beyond that point, while remaining piecewise linear and therefore analytically tractable.

%%%%%%%%%%%%%%%%%%%%%%%

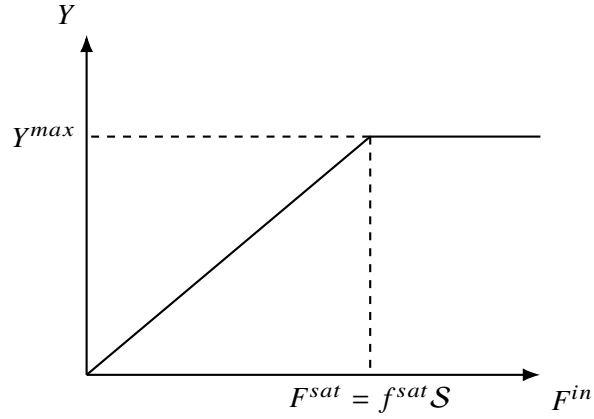
\begin{figure}[tbp]
\centering
\begin{tikzpicture}[
  scale=1.5,
  >=Latex,
  axis/.style={->, thick},
  curve/.style={thick},
  guide/.style={dashed, thick}
]
  % Axes
  \draw[axis] (0,0) -- (4,0) node[below right] {$F^{in}$};
  \draw[axis] (0,0) -- (0,3) node[above left] {$Y$};

  % Key points
  \coordinate (O) at (0,0);
  \coordinate (S) at (2.5,2.1);
  \coordinate (R) at (4,2.1);

  % Curve
  \draw[curve] (O) -- (S) -- (R);

  % Guides
  \draw[guide] (S) -- (2.5,0);
  \draw[guide] (S) -- (0,2.1);

  % Labels
  \node[below] at (2.5,0) {$F^{sat}=f^{sat}{\cal S}$};
  \node[left]  at (0,2.1) {$Y^{max}$};
\end{tikzpicture}
\caption{Piecewise-linear saturating production function used to model crop and grassland yields as a function of total nitrogen input $F^{in}$. Production increases linearly in the unsaturated regime up to the saturation threshold $F^{sat}=f^{sat}{\cal S}$, then reaches a plateau at the maximum yield $Y^{max}$.}
\label{fig:yield_function}
\end{figure}

%%%%%%%%%%%%%%%%%%%%%%%%ù

%\begin{center}
%Piecewise-linear saturating production function used to model crop and grassland yields as a function of total nitrogen input $F^{in}$. Production increases linearly in the unsaturated regime up to the saturation threshold $F^{sat}=f^{sat}{\cal S}$, then reaches a plateau at the maximum yield $Y^{max}$.
%\end{center}
\medskip

\noindent\textbf{Unsaturated regime.}
Let
\[
F^{in}_{CROP} := F^{synth}_{CROP} + BNF_{CROP} + F_{EXC \to CROP},
\]
\[
F^{in}_{GRASS} := F^{synth}_{GRASS} + BNF_{GRASS} + F_{EXC \to GRASS}.
\]
denote the total nitrogen inputs to cropland and grassland, respectively.
If
\[
F^{in}_{CROP} < f^{sat}_{CROP}\,{\cal S}_{CROP},
\]
crop production is given by
\begin{equation}
Y_{CROP} := F_{CROP \to CRP} = \gamma_{CROP}\, F^{in}_{CROP}.
\end{equation}
Similarly, if
\[
F^{in}_{GRASS} < f^{sat}_{GRASS}\,{\cal S}_{GRASS},
\]
grassland production satisfies
\begin{equation}
Y_{GRASS} := F_{GRASS \to GRP} = \gamma_{GRASS}\, F^{in}_{GRASS}.
\end{equation}
The unsaturated yield coefficients satisfy $0<\gamma_{CROP}, \gamma_{GRASS} \le 1$.
At stationarity, nitrogen balance then implies
\[
F_{CROP \to \emptyset} = (1-\gamma_{CROP})\,F^{in}_{CROP}, \qquad
F_{GRASS \to \emptyset} = (1-\gamma_{GRASS})\,F^{in}_{GRASS}.
\]

\noindent\textbf{Saturated regime.}
If
\[
F^{in}_{CROP} \ge f^{sat}_{CROP}\,{\cal S}_{CROP},
\]
crop production reaches its maximum value
\begin{equation}
F_{CROP \to CRP} = Y^{max}_{CROP} := \gamma_{CROP}\, f^{sat}_{CROP}\,{\cal S}_{CROP}.
\end{equation}
Likewise, if
\[
F^{in}_{GRASS} \ge f^{sat}_{GRASS}\,{\cal S}_{GRASS},
\]
then
\begin{equation}
F_{GRASS \to GRP} = Y^{max}_{GRASS} := \gamma_{GRASS}\, f^{sat}_{GRASS}\,{\cal S}_{GRASS}.
\end{equation}
Excess nitrogen is routed to environment ($\emptyset$):
\[
F_{CROP \to \emptyset} = F^{in}_{CROP} - Y^{max}_{CROP}, \qquad
F_{GRASS \to \emptyset} = F^{in}_{GRASS} - Y^{max}_{GRASS}.
\]

%%%%%%%%%%%%%%%%%%%%%%%%%%%%
\subsection{Constraints}  \label{subsection:constraints}
%%%%%%%%%%%%%%%%%%%%%%%%%%%%

There are many different ways to constrain the system. We may roughly distinguish three types of constraints:

-- {\em demand stock requirements}.  Human population and livestock stocks, denoted $v_{POP}$ and $v_{LVS}$, are treated as stationary quantities. Fixing both quantities exactly is generally incompatible with the linear stationarity conditions, unless synthetic inputs and trade flows are adjusted accordingly;

-- {\em reverse allocation ratios}, mainly human and livestock diets, for which one may. 
 prescribe reference values,
\[
\begin{aligned}
\left(
\frac{F_{CRP \to POP}}{v_{POP}},
\frac{F_{LVP \to POP}}{v_{POP}}
\right)
&\doteq (\alpha_{POP}, 1-\alpha_{POP}), \\
\left(
\frac{F_{CRP \to LVS}}{v_{LVS}},
\frac{F_{GRP \to LVS}}{v_{LVS}}
\right)
&\doteq (\alpha_{LVS}, 1-\alpha_{LVS}).
\end{aligned}
\]

Technical parameter values representative of contemporary 2014 French agriculture are provided in Suppl. Mat.

%%%%%%%%%%%%%%%

\subsection{Numerical resolution of flux constraints}  \label{subsection:num-resol}

%%%%%%%%%%%%%%%%

In this work, two methods were used  to solve the constraints. The first one is based on a  resolution of the flux constraints by matrix inversion; see Suppl. Mat. The 
matrix represents flux conservation at each node, as well as demand stock requirements and allocation
constraints (\S \ref{subsection:constraints}). 

The heterogeneous nature of its coefficients, and their mixed signs, suggest no specific method of resolution besides direct {\em numerical inversion}.

%%%%%%%%ù

\subsection{Hierarchical formulas}  \label{subsection:HF}

%%%%%%%%%%%%%%%ù

By contrast, the second method allows a  {\em formal approximate resolution} in terms of stock variables. By formal, we mean
that we obtain the solution as a function $f({\bf x})$ of the set ${\bf x} = \{\textsc{VolRate}$, 
$\textsc{CollRate}$, ${\cal S}_{CROP}$, ${\cal S}_{GRASS}$, $k_{CRP\to LVS}$, $k_{EXC\to CROP}/k_{EXC\to GRASS}\}$,
 of main variables of the problem, see \S \ref{subsection:main-variables}. 
The function takes  the form of compact, smooth {\bf hierarchical formulas}, which are  good approximations of the numerical solution  inside geometric domains $\cal D$ in ${\bf x}$ defined by sets of inequalities. Hierarchical formulas    are usually discontinuous at the boundary ${\cal D}\cap{\cal D}'$ between two domains ${\cal D},{\cal D}'$, where the solution changes rapidly. We interpret this by saying that geometric domains 
partition the parameter set $\bf x$ into distinct {\bf metabolic phases}. 

\begin{figure}[tbp]
    \centering
    \includegraphics[width=\textwidth]{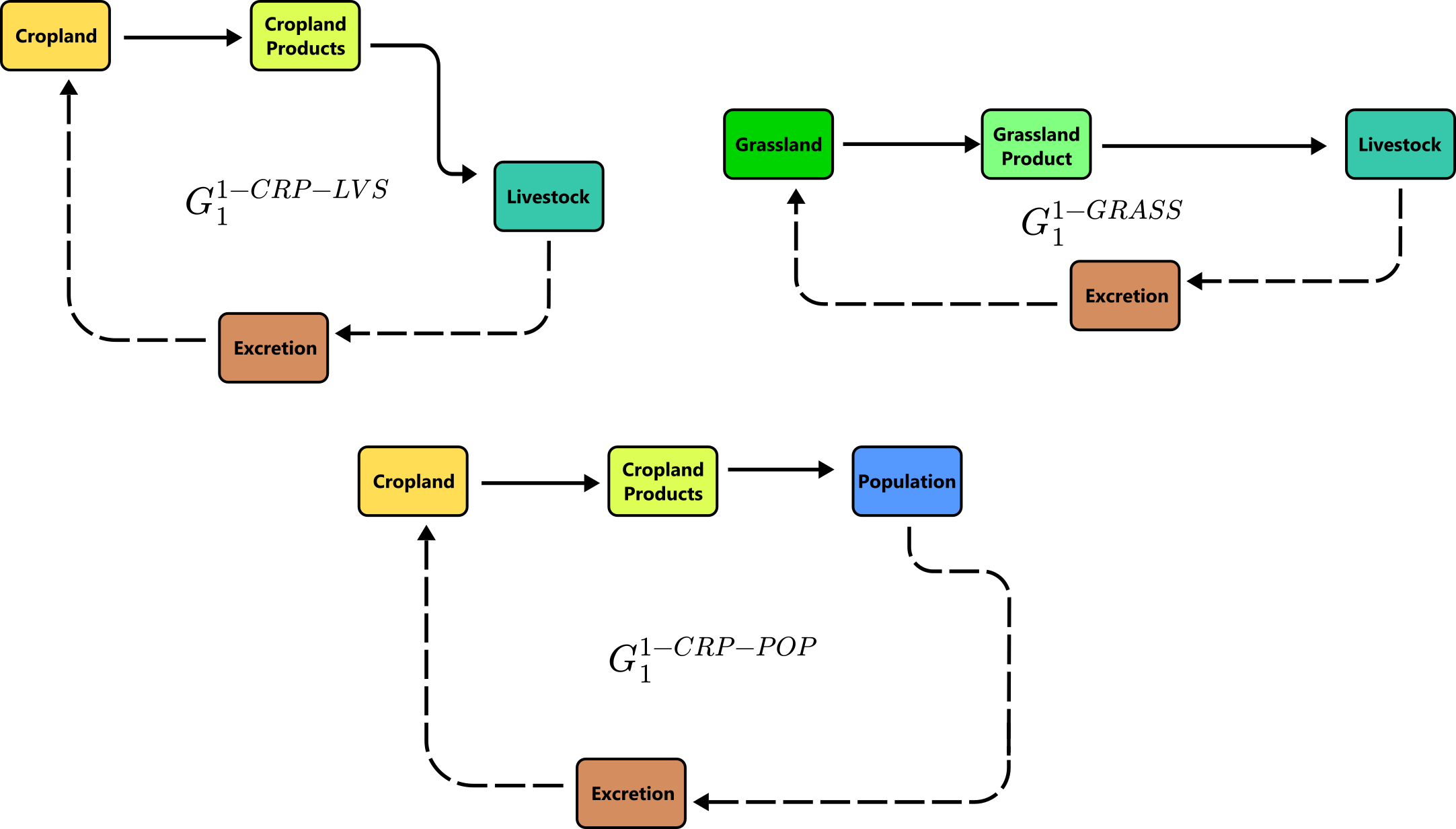}
    \caption{The three simplest possible operating cycles}
    \label{fig:three-simplest-possible-cycles}
\end{figure}

\medskip All these domains, save one, are characterized by the presence of one or several {\bf operating metabolic cycles}. Simplest  metabolic cycles are cycles in the usual sense (topological cycles); there are three instances of such cycles, see Fig. \ref{fig:three-simplest-possible-cycles},

\begin{itemize}
    \item[\textbullet] 
 the {\bf crop-livestock driven metabolic cycle} $\{CROP \to CRP \to 
LVS \to EXC \to CROP\}$ with the return edge $EXC\to CROP$ (denoted $G_1^{1-CRP-LVS}$). In phase $\mathsf{CE.1}$,  %{\color{red}SS.1.3}, 
$G_1^{1-CRP-LVS}$ is {\bf operating}, which means that alternative edges leaving the cycle, namely, 
$CROP\to \emptyset$,  $CRP\to POP$,  $EXC\to \emptyset$ and $EXC\to GRASS$, have a small direct allocation ratio attached to them. The
condition  for the outflow $CROP\to \emptyset$ always holds true, but the other ones imply inequalities for $\bf x$, 
namely, 
\begin{equation} 
{\textsc{VolRate}}<0.5,\qquad  k_{CRP\to LVS}>0.5, \qquad k_{EXC\to GRASS} < k_{EXC\to CROP}
\end{equation}
We say that allocation $EXC\to CROP$ of excreta to CROP  is {\bf dominant}, and so is the allocation $CRP\to LVS$ of crops to livestock;

 \item[\textbullet]  when   ${\textsc{VolRate}}<0.5$, and  $EXC\to GRASS$
is dominant (defining phase $\mathsf{GE}$) , the {\bf grass driven metabolic cycle} $G_1^{1-GRASS} = \{GRASS \to GRP \to LVS\to EXC\to GRASS\}$ is operating.  

 \item[\textbullet]  when   ${\textsc{VolRate}}<0.5$,  ${\textsc{CollRate}}>0.5$, and  $CRP\to POP$,  $EXC\to CROP$
are dominant (defining phase $\mathsf{CE.3}$) , the {\bf crop product-population driven metabolic cycle} $G_1^{1-CRP-POP} = \{CROP \to CRP \to POP\to EXC\to CROP\}$ is operating.

\end{itemize}

There also appear more complex metabolic cycles, made up of intertwined topological cycles. 

\medskip When operating, metabolic cycles $G$ have interesting properties:
\begin{itemize}
    \item[(1)] {\em stocks of $G$-compartments are to leading order independent of external rates} 
    connecting compartments which are not in $G$;
    \item[(2)] the {\em relative proportion
of stocks of  $G$-compartments is to leading order independent of the rates of edges entering or leaving $G$}, 
and may be computed assuming that these vanish;
\item[(3)] the N circulation within $G$ {\em maintains itself without
any external influx for a large period of time};
\item[(4)] there is a {\em large multiplier effect} attached to any influx into $G$ (due e.g. to supplementary
fertilization), namely, varying an influx by 1 unit typically leads to supplementary stocks of $G$-compartments of 
order $Z_G^{-1}$, where $Z_G^{-1}\gg 1$ (called: {\bf $G$-weight}) is a large constant depending both on the structure of $G$ and 
on the way that $G$ sits inside the whole network (but not on external rates). 
\end{itemize}

Such properties point at a {\em nearly autonomous functioning} of the compartments of $G$ as a whole. If emphasis is laid on {\em core compartments}, i.e. compartments of $G$, then property (1) may be
referred to as {\bf robustness}: {\em core compartments are unaffected by external fluxes}. Property (4), on the other hand, is a type of {\bf resilience}: {\em any not-too-large  perturbation within $G£$ may be taken care of by a small influx variation}. The sensitivity to rates of edges leaving $G$ ({\bf parameter sensitivity}) may be assessed by looking at the formula for $Z_G$. 

\medskip On the contrary,
these properties are wrong when $G$ is {\bf broken}, i.e. non operating.
For instance,  large multiplier effects (4) are absent from the {\em trivial phase}, in which all metabolic cycles are
broken. {\em Any type of perturbation breaking a metabolic cycle} (or, possibly, substituting another operating metabolic
cycle) {\em causes a large-scale rearrangement of the network}, therefore  leading-order shifts in stocks and fluxes;
the network is not resilient to such perturbations.

\medskip Roughly speaking, in the presence of a unique operating metabolic cycle $G$, the main characteristics of the 
system may be obtained from the relative proportion of stocks inside $G$, the $G$-weight $Z_G^{-1}$, and an {\bf 
effective graph} obtained by summing the stocks of all $G$-compartments, replacing them by a unique compound
stock indexed by $G$, and redefining according to simple rules the kinetic rates to and from $G$. The graphical procedure is called  {\bf merging}, and the rate redefinition a {\bf renormalization step}. 

\medskip The {\bf hierarchical algorithm} presented in Section 3 of Suppl. Mat. discusses these rules under very general assumptions,
allowing the more involved analysis of cases when multiple metabolic cycles compete. Several successive renormalization steps are needed in the GRAFS model, depending on the phase, producing a {\bf hierarchical, multi-level structuration}
of the network in general. Relevant quantities, like stationary stocks, are then approximated through recursive
 hierarchical formulas. The above properties (1--4) and associated concepts of robustness and resilience, and parameter sensitivity, are discussed in  particular in \S 3.5 of 
 Suppl. Mat.

\medskip Compared to direct numerical resolution (\S \ref{subsection:num-resol}), this algorithm has the clear advantage of yielding 
interpretable results which can be directly  used to inform decision making. Properties (1)-(4), and  the concepts sketched thereafter, are discussed at a more formal and quantitative level in Suppl. Mat. The obvious drawback is that, while exact in a certain abstract limit, it is at best semi-quantitative when using 'real' parameter values, with relative errors in practice as high as $30\%$. Yet general tendencies - as confirmed by comparison with numerics - are {\em always} reliable (see Suppl. Mat., Section 8). In this sense, it can be viewed as a decision tool for N metabolism, in combination with numerical tools.

%%%%%%%%%%%%%%%%%%%%%%%%%%%
%%%%%%%%%%%%%%%%%%%%%%%%%%

\section{Results and discussion}
\label{section:results}
%%%%%%%%%%%%%%%%%%%%%%%%%%
%%%%%%%%%%%%%%%%%%%%%%%

\subsection{Classification of metabolic phases}

The analytical exploration of the network reveals the existence of several distinct
\emph{metabolic phases}, each corresponding to a specific operating mode of the
agro-food system.
These phases differ  in nitrogen circulation
patterns, allocation priorities, and dominant recycling pathways, and are characterized by an organization and 
sensitivity properties discussed sketchily in \S \ref{subsection:HF} and in more detail in Suppl. Mat.

\medskip We distinguish several broad classes of metabolic phases, that we call {\em regimes}:

\begin{itemize}
\item[(i)] \emph{Low-volatilization} regime, characterized by
\[
\textsc{VolRate}
<
k_{EXC \to CROP} + k_{EXC \to GRASS},
\]
as opposed to \emph{high-volatilization} regime, where
\[
\textsc{VolRate}
>
k_{EXC \to CROP} + k_{EXC \to GRASS}.
\]    Equivalently, \textsc{VolRate}$\gtrless 0.5$;

\item[(ii)] \emph{population-driven} regime, defined by
\[
k_{CRP \to POP} > k_{CRP \to LVS},
\]
as opposed to \emph{livestock-driven} regime, for which
\[
k_{CRP \to LVS} > k_{CRP \to POP},
\]
with intermediate configurations corresponding to approximately even crop allocation.

\item[(iii)] \emph{cropland-driven} regime, for which
\[
k_{EXC \to CROP} > k_{EXC \to GRASS},
\]
as opposed to \emph{grassland-driven} regime, where
\[
k_{EXC \to GRASS} > k_{EXC \to CROP}.
\]

\item[(iv)] the \emph{high-recycling} regime, defined by
\[
\textsc{CollRate}
>
 0.5
\]
and the \emph{low-recycling} regime,
\[
\textsc{CollRate}
<
0.5
\]

\end{itemize}

%\begin{equation} \label{Fig:Decision_tree}
%\includegraphics[scale=0.8]{Decision_tree}   
%\end{equation}
%\begin{center}
%{\textsc{Fig. \ref{Fig:Decision_tree} -- Decision tree.}}
%\end{center}

\begin{figure}[tbp]
    \centering
    \includegraphics[width=0.6\textwidth]{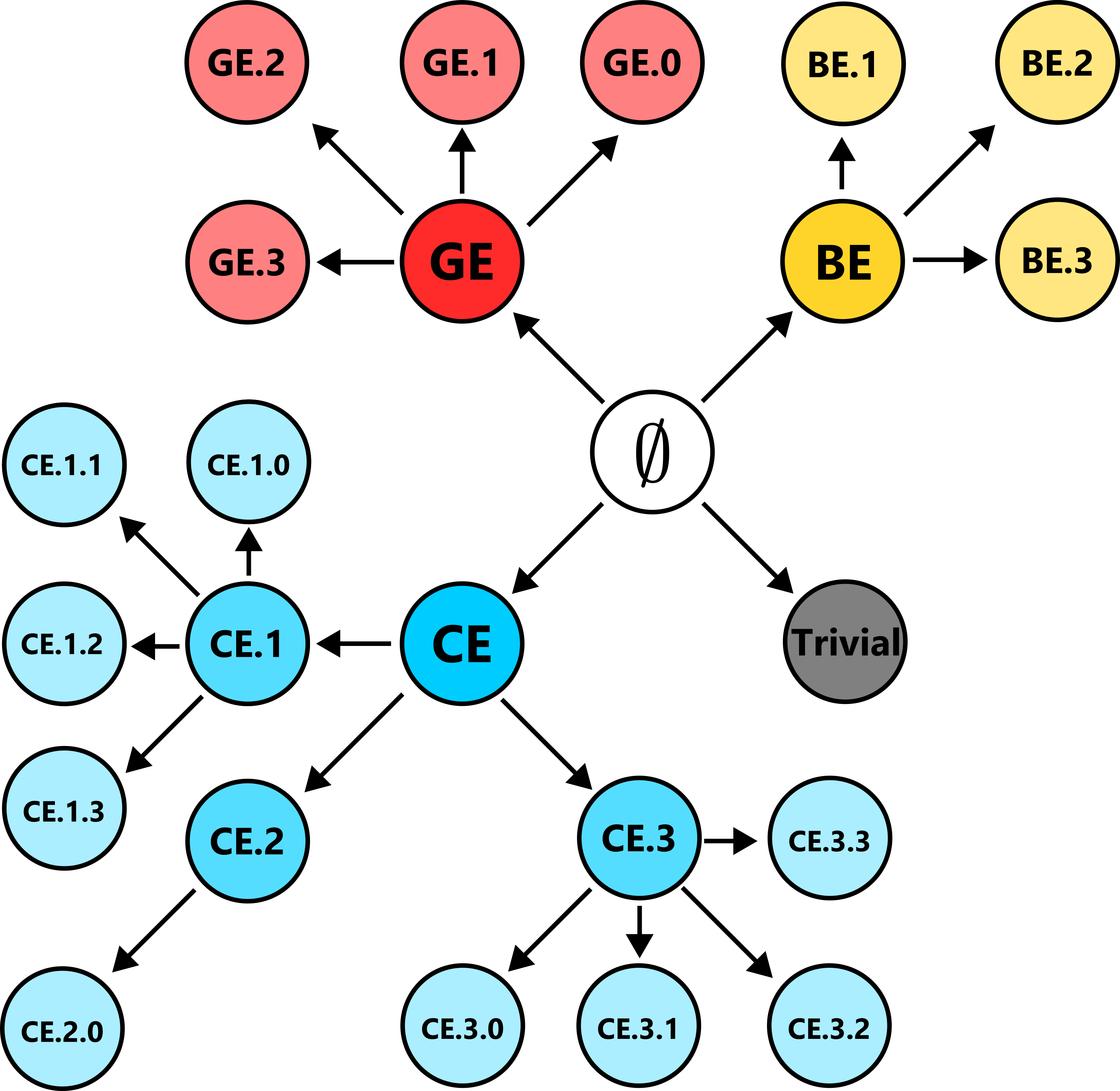}
    \caption{Analytical classification tree of the metabolic phases of the agro-food nitrogen network. Nodes denote qualitatively distinct operating modes defined by excreta allocation, crop allocation, and the balance between recycling and losses; the Trivial node corresponds to weakly structured configurations.}
    %\caption{Mapping of the metabolic phases of the agro-food nitrogen network. Each node represents a qualitatively distinct operating mode, defined by dominant allocation pathways and loss or recycling regimes. The node \emph{Trivial} corresponds to weakly structured configurations dominated by high volatilization or low recycling. Non-trivial phases are obtained through a hierarchical sequence of decision rules based on (i) excreta allocation between cropland and grassland, (ii) allocation of cropland products between human consumption and livestock feeding, and (iii) the balance between recycling and losses. Colors indicate major phase families (BE: balanced excreta allocation; CE and GE: cropland- or grassland-driven excreta allocation, respectively), while sublabels denote refined subclasses arising from secondary allocation criteria. The tree provides a compact representation of the analytically derived phase structure of the metabolic network.}
    \label{Fig:Decision-tree}
\end{figure}

\paragraph{Decision rules.}  In many cases, operating metabolic cycles depend on the choice of a dominant allocation edge for a compartment. When two competing allocation rates are comparable, the cycle contains both. This transitional regime is singled out using 
 a dimensionless decision threshold
\BEQ
\texttt{Decision\_thr} = 1.3
\EEQ
for the ratios.

For instance, when comparing the excreta allocation toward cropland and grassland,
we classify the outcome as follows:
\begin{itemize}
\item[\textbullet] $k_{EXC \to CROP} \succ k_{EXC \to GRASS}$ if
$\displaystyle \frac{k_{EXC \to CROP}}{k_{EXC \to GRASS}} > \texttt{Decision\_thr}$;
\item[\textbullet] $k_{EXC \to CROP} \prec k_{EXC \to GRASS}$ if
$\displaystyle \frac{k_{EXC \to GRASS}}{k_{EXC \to CROP}} > \texttt{Decision\_thr}$;
\item[\textbullet] $k_{EXC \to CROP} \sim k_{EXC \to GRASS}$ otherwise.
\end{itemize}

\paragraph{Trivial phase.}
We first identify an extended \emph{trivial phase}, characterized by weak internal
structuring of nitrogen flows and a multifaceted boundary in parameter space.
This phase is selected whenever one of the following conditions is satisfied:
\begin{itemize}
\item[\textbullet] $\textsc{VolRate} > 0.5$ (high-volatilization regime);
\item[\textbullet] $\textsc{VolRate} < 0.5$,
$k_{CRP \to LVS} \preceq k_{CRP \to POP}$, 
$k_{EXC \to CROP} \succ k_{EXC \to GRASS}$, $\textsc{CollRate} < 0.5$ 
(low-volatilization, human-driven or even crop allocation, cropland-driven, low-recycling regime)
\end{itemize}

It is interesting to note at this stage that :

\begin{itemize}
    \item[\textbullet]
  the volatilization rate is the key variable; current French agriculture has a relatively low
$\textsc{VolRate} = 0.32 < 0.5$;

\item[\textbullet]
being in a low-recycling regime  ($\textsc{CollRate} < 0.5$ in most contemporary agricultural systems) does not automatically imply  that the network is in its trivial phase; this holds only if human vegetal nutrition  (the network edge $CRP\to POP$) is part of any dominant metabolic cycle, which excludes  the case of the  livestock-driven or grassland-driven regimes (in which one or both of the cycles $G_1^{1-CRP-LVS}$ or $G_1^{1-GRASS}$ are dominant, see Suppl. Mat. for more details).

\end{itemize}

\medskip
\noindent
\emph{If the system does not fall within the trivial phase}, the following hierarchical
decision rules are applied to determine the metabolic phase (see Fig. \ref{Fig:Decision-tree}, and Suppl. Mat., Fig. 4-6 for a detailed illustration including operational metabolic cycles):

\medskip
\noindent\textbullet\ \textbf{Primary split (excreta allocation)}
\begin{itemize}
\item[\textbullet] $\emptyset \to BE$ {\em (Balanced Excretion)} if $k_{EXC \to CROP} \sim k_{EXC \to GRASS}$;
\item[\textbullet] $\emptyset \to CE$ {\em (Cropland Excretion)} if $k_{EXC \to CROP} \succ k_{EXC \to GRASS}$;
\item[\textbullet] $\emptyset \to GE$ {\em (Grassland Excretion)} if $k_{EXC \to CROP} \prec k_{EXC \to GRASS}$.
\end{itemize}

\medskip
\noindent\textbullet\ \textbf{Secondary split within BE (crop product allocation)}
\begin{itemize}
\item[\textbullet] $BE \to BE.1$ if $k_{CRP \to LVS} \succ k_{CRP \to POP}$;
\item[\textbullet] $BE \to BE.2$ if $k_{CRP \to POP} \sim k_{CRP \to LVS}$;
\item[\textbullet] $BE \to BE.3$ if $k_{CRP \to LVS} \prec k_{CRP \to POP}$.
\end{itemize}

\medskip
\noindent\textbullet\ \textbf{Secondary split within CE (crop product allocation)}
\begin{itemize}
\item[\textbullet] $CE \to CE.1$ if $k_{CRP \to LVS} \succ k_{CRP \to POP}$;
\item[\textbullet] $CE \to CE.2$ if $k_{CRP \to POP} \sim k_{CRP \to LVS}$;
\item[\textbullet] $CE \to CE.3$ if $k_{CRP \to LVS} \prec k_{CRP \to POP}$.
\end{itemize}

\medskip
\noindent\textbullet\ \textbf{Secondary split within GE (loss-dominated vs.\ recycling-dominated)}
\begin{itemize}
\item[\textbullet] $GE \to GE.0$ if $k_{EXC \to CROP} + k_{LVS\to LVP} <
k_{GRASS \to \emptyset} + \textsc{VolRate};$ 
\item[\textbullet] $GE \to GE\setminus GE.0$ else.
\end{itemize}

\medskip
\noindent\textbullet\ \textbf{Secondary split within $\mathsf{GE}\setminus \mathsf{GE}.0$ (crop product allocation)}
\begin{itemize}
\item[\textbullet] $GE\setminus GE.0 \to GE.1$ if $k_{CRP \to LVS} \succ k_{CRP \to POP}$;
\item[\textbullet] $GE\setminus GE.0 \to GE.2$ if $k_{CRP \to POP} \sim k_{CRP \to LVS}$;
\item[\textbullet] $GE\setminus GE.0 \to GE.3$ if $k_{CRP \to LVS} \prec k_{CRP \to POP}$.

\end{itemize}

\medskip \noindent For France 2014, we have $k_{EXC \to CROP} = 0.52$, $k_{EXC \to GRASS}=0.27$, $\textsc{CollRate}=0.15$ and $k_{CRP \to LVS} \succ k_{CRP \to POP}$. This means that France 2014 is in  phase $\mathsf{CE.1}$.

\begin{figure}[tbp]
\centering
\includegraphics[width=\textwidth]{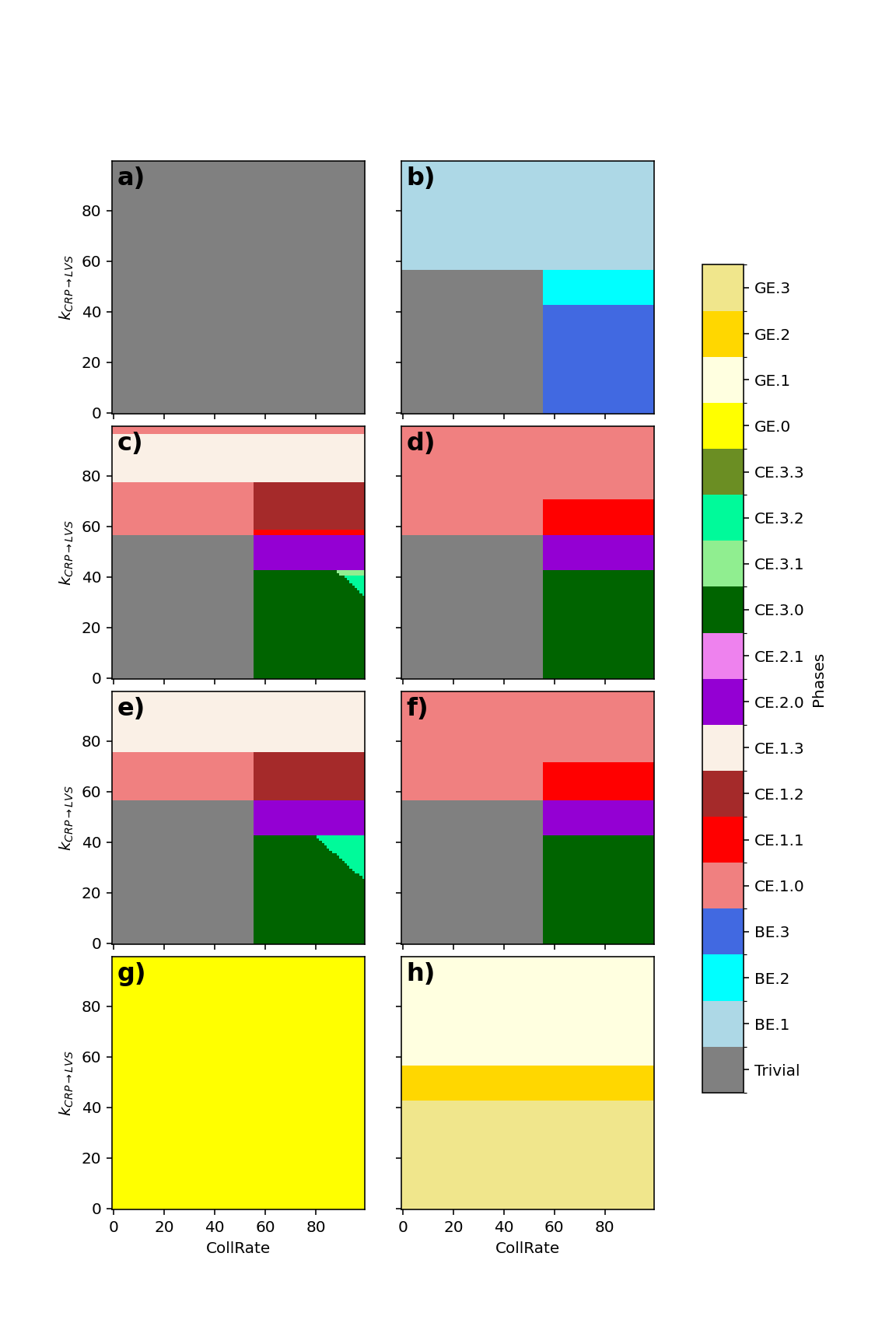}
\caption{Phase diagrams in the space of excreta collection and  crop allocation  for eight combinations of excreta volatilization and excreta allocation: (a) (0.52, 0.50), (b) (0.32, 0.50), (c) (0.20, 0.60), (d) (0.20, 0.90), (e) (0.15, 0.60), (f) (0.15, 0.95), (g) (0.42, 0.30), and (h) (0.15, 0.30), where each pair is given as ($\textsc{VolRate}$, $\textsc{ExcAllocRatio}$). High volatilization collapses the system toward the trivial regime, whereas lower volatilization and more polarized recycling reveal differentiated balanced- ($\mathsf{BE}$), cropland- 
($\mathsf{CE}$), and grassland-excretion ($\mathsf{GE}$) phases.}
\label{fig:phase-interaction}
\end{figure}

\medskip
\noindent\textbf{Phase diagrams.}
Phase diagrams are used to visualize the classification of metabolic phases across the
parameter space defined by the collecting rate \textsc{CollRate}   (horizontal axis) and crop allocation ratio $k_{CRP\to LVS}$, (vertical axis), for fixed values of volatilization
rate \textsc{VolRate} and excreta allocation ratio
$$\textsc{ExcAllocRatio}=\frac{k_{EXC\to CROP}}{k_{EXC\to CROP}+k_{EXC\to GRASS}}$$
 See Fig.
\ref{fig:phase-interaction}.

\subsection{French agriculture, three case studies}

%%%%%%%%%%%%%%%%%%%%%%%%%%%
%%%%%%%%%%%%%%%%%%%%%%%%%%%%%%%

\Bigskip We now illustrate the framework with three case studies for French agriculture. Each case study explores a different set of control variables and constraints, using the numerical resolution procedure described in Suppl. Mat.  Synthetic inputs to grassland are scaled proportionally to cropland inputs using the French reference ratio $f^{synth}_{GRASS}/f^{synth}_{CROP}=0.04$.

%%%%%%%%%%%%%%%%%%%%%%%%%%%

\subsubsection{A - Substitutability between synthetic fertilization and biological nitrogen fixation}
\label{subsection:natural-synthetic}

We examine to what extent biological nitrogen fixation (BNF) by legumes can substitute for synthetic nitrogen fertilization while keeping the French land base, livestock system, and diets fixed. Cropland and grassland areas are held constant at their French reference values, ${\cal S}_{CROP}=18.3$~Mha and ${\cal S}_{GRASS}=9.1$~Mha. The livestock population is fixed at $LU=20.25\times 10^6$, and the human stock is fixed to the French 2014 population, rounded to $67\times 10^6$ inhabitants. Human and livestock diets (see (\ref{eq:diet})) are kept at their French reference composition, as are all technical coefficients, including excreta collection and recycling parameters (see Suppl. Math., Section 1). Under this closure, boundary fluxes of cropland products, grassland products, and livestock products are left free.

The two control variables are (i) the mean synthetic nitrogen input to cropland, $f^{synth}_{CROP}$ (kgN\,ha$^{-1}$), and (ii) the share of legumes in cropland, $\tau_{LEG,CROP}$ (\% of cropland area). Cropland BNF per ha is assumed to increase linearly with legume share as $C_{T,CROP}~\tau_{LEG,CROP}$, with $C_{T,CROP}=264$~kgN\,ha$^{-1}$, while grassland BNF is kept constant at $C_{T,GRASS}=38.9$~kgN\,ha$^{-1}$. 
Because cropland products are used both for human food and livestock feed, feeding capacity is not computed directly from the total net cropland-product boundary flux. Domestic cropland production is first allocated to the fixed livestock herd. The remaining domestic cropland products are then converted into a human-equivalent feeding capacity:
\[
B^{avail}_{POP}
=
\max\left(0,\,
F_{CROP\rightarrow CRP}
-
F_{CRP\rightarrow LVS}
\right),
\]
and
\[
\rho_{POP}
=
\frac{B^{avail}_{POP}}{{\mathrm{conv}}_{POP}\,{\cal S}_{CROP}},
\]
where conv$_{POP}=6.44$~kgN\,cap$^{-1}$\,yr$^{-1}$ is the total annual nitrogen requirement per inhabitant. The indicator is therefore a dietary N-equivalent feeding-capacity density based on domestic cropland products remaining after feeding the fixed livestock herd.

\begin{figure}[tbp]
\centering
\includegraphics[width=0.9\textwidth]{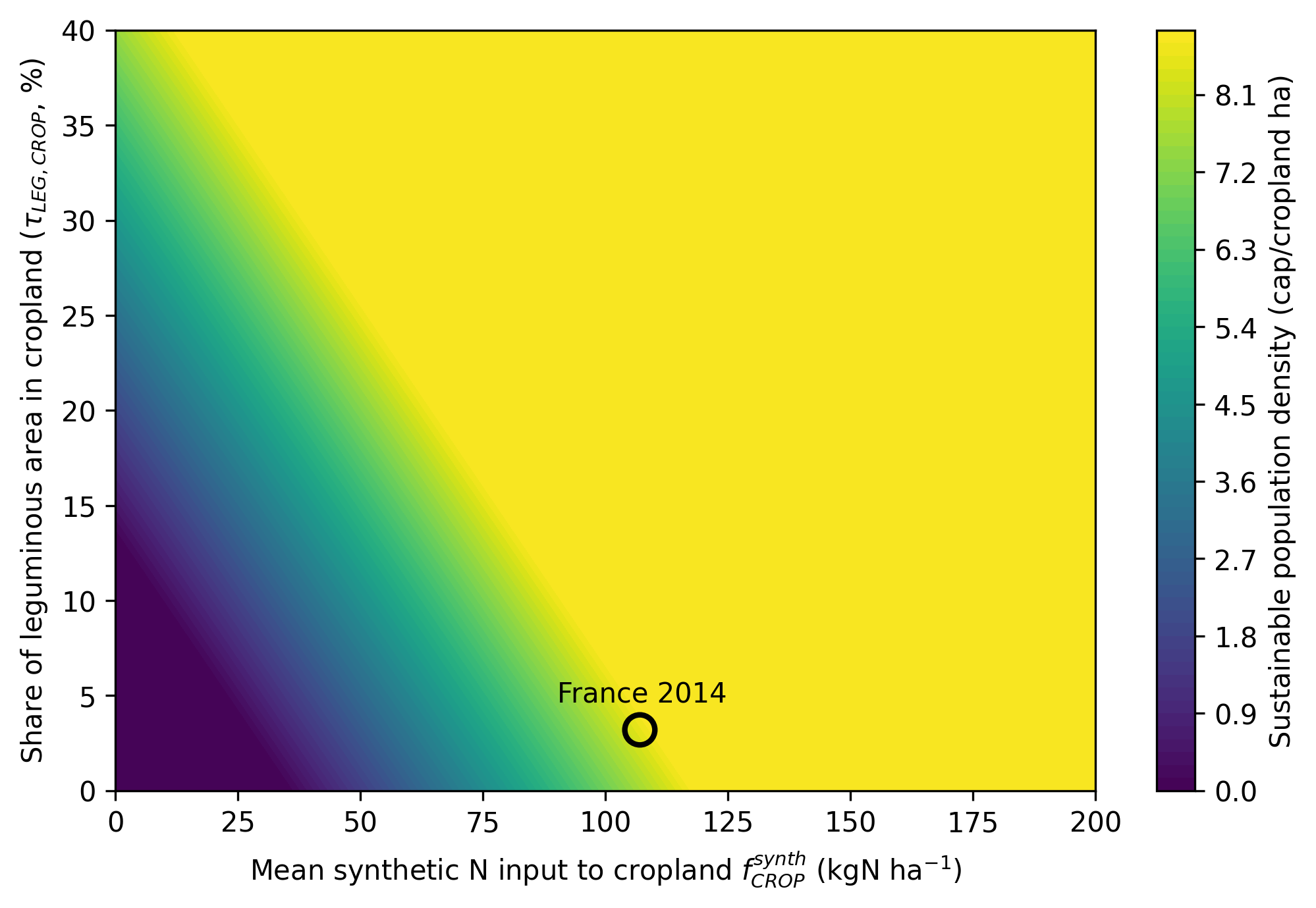}
\caption{Dietary N-equivalent feeding-capacity density $\rho_{POP}$ (capita\,ha$^{-1}$ cropland) as a function of mean synthetic N input to cropland and legume share in cropland. Cropland and grassland areas, herd size, French 2014 population, human and livestock diets, and technical coefficients are fixed. Boundary fluxes of cropland products, grassland products, and livestock products are free. Feeding capacity is computed from domestic cropland products remaining after feeding the fixed livestock herd and converted using the total annual human nitrogen requirement. The circle marks the French 2014 reference configuration. The whole explored domain belongs to phase $\mathsf{CE.1.0}$.}
\label{fig:1st_case}
\end{figure}

The resulting heatmap (Figure~\ref{fig:1st_case}) shows that $\rho_{POP}$ increases with both synthetic fertilization and legume-based fixation before reaching a broad high-input plateau. Iso-capacity contours are diagonal in the low-input domain, reflecting the direct substitutability between synthetic nitrogen and BNF in the effective cropland nitrogen input. Increasing the legume share can therefore compensate for lower synthetic fertilization, up to the point where cropland production reaches its saturation threshold.

At the French 2014 reference point, defined by $f^{synth}_{CROP}=107.0$~kgN\,ha$^{-1}$ and $\tau_{LEG,CROP}=3.2\%$, the model yields $\rho_{POP}=8.81$~cap\,ha$^{-1}$ cropland, corresponding to a total feeding capacity of 161.27 million inhabitants in dietary N-equivalent terms. This value is higher than the fixed French population because, at the reference input level, domestic cropland production exceeds the cropland-product requirement of the fixed livestock herd and leaves an exportable surplus that can be expressed as additional human-equivalent food supply.

The heatmap also illustrates the substitutability between synthetic fertilization and legume-based BNF. Starting from the French reference point, reducing synthetic fertilization by 20\% would require increasing the legume share in cropland to about 11\% in order to maintain the same feeding-capacity density. This corresponds to approximately 2.0~Mha of cropland under legumes, compared with about 0.6~Mha at the current 3.2\% legume share.

%In low-input configurations, domestic cropland production may be insufficient to cover the cropland-product demand associated with the fixed population and livestock herd. The balance is then closed by net imports of cropland products. These imports maintain the prescribed stocks in the stationary solution, but they are not credited in the feeding-capacity indicator, which measures only the domestic cropland products available for human-equivalent food supply once the fixed crop-feed requirement of livestock has been accounted for. As a result, the indicator tends toward zero in very low-production configurations rather than becoming negative.

The whole explored domain is classified as $\mathsf{CE.1.0}$. The prefix $\mathsf{CE}$ indicates a crop-excretion regime, in which recycled excreta are preferentially routed toward cropland rather than grassland. The branch $\mathsf{CE.1}$ corresponds to a livestock-oriented allocation of cropland products, characterized by the crop-livestock driven metabolic cycle $G_1^{1-CRP-LVS}$, and consistent with the fixed herd and the large crop-feed requirement in the French reference system. Loss rates in the system, in particular, the low value of \textsc{CollRate}, prohibit a multi-level metabolic structuring, and account for the further $\mathsf{.0}$ in the domain name. 

%%%%%%%%%%%%%%%%%%%%%%%%%%%

\subsubsection{B - Joint effects of livestock population and synthetic fertilization on sustainable population density}
\label{subsection:LU-Fsynth}

Case study~B uses the same closure and feeding-capacity indicator as Case study~A, but varies livestock population instead of legume cover. The legume share is fixed at its French reference value, while the two control variables are the livestock population $LU$ and the mean synthetic nitrogen input to cropland, $f^{synth}_{CROP}$. For each parameter pair, the model keeps the French 2014 population, diets, land areas and technical coefficients fixed, and lets the three product border fluxes $\Phi_{CRP}$, $\Phi_{GRP}$ and $\Phi_{LVP}$ adjust endogenously.

\begin{figure}[tbp]
\centering
\includegraphics[width=0.9\textwidth]{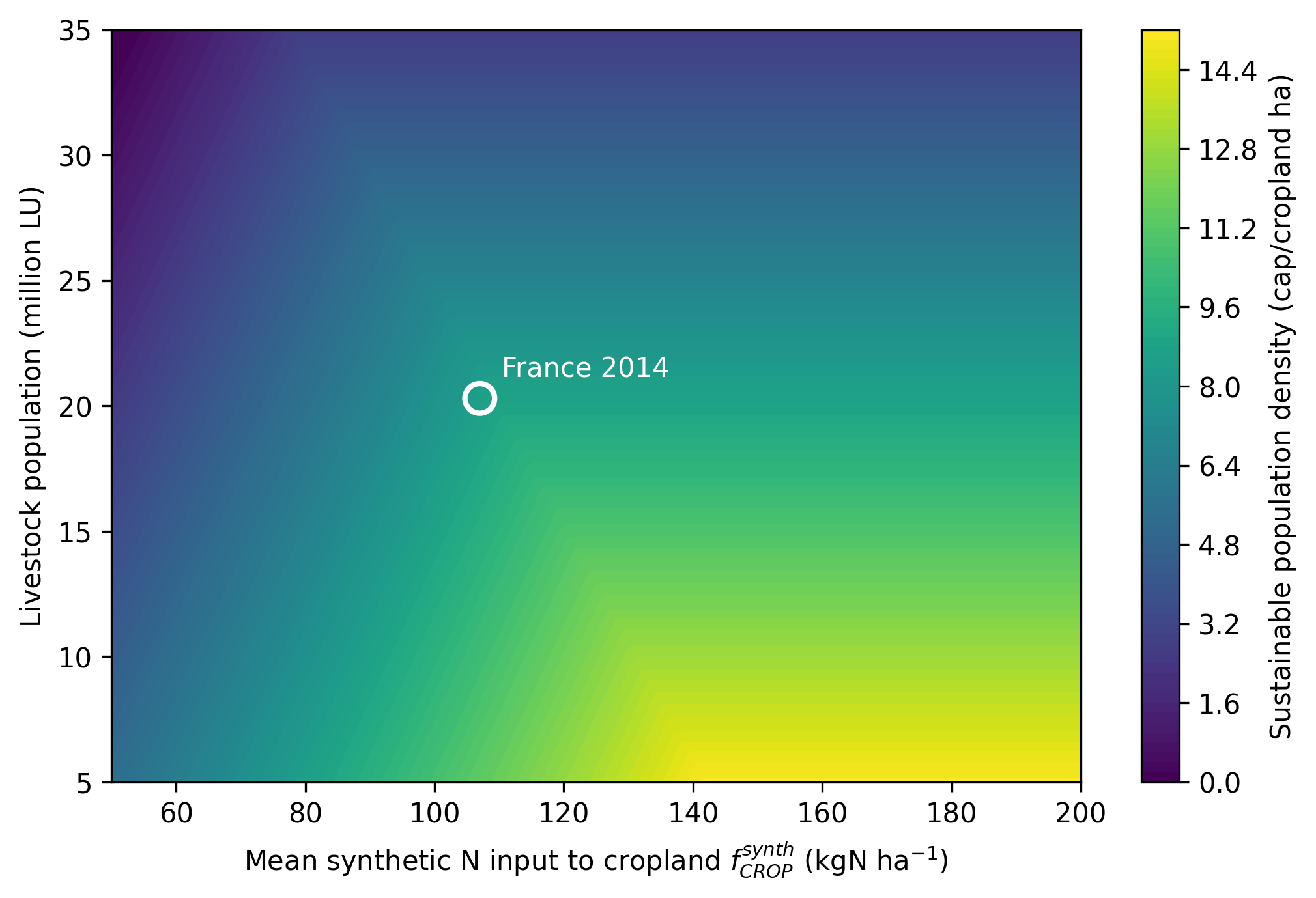}
\caption{Dietary N-equivalent feeding-capacity density $\rho_{POP}$ (capita\,ha$^{-1}$ cropland) as a function of livestock population and mean synthetic N input to cropland. Cropland and grassland areas, legume share, French 2014 population, human and livestock diets, and technical coefficients are fixed. Boundary fluxes of cropland products, grassland products, and livestock products are free. Feeding capacity is computed from domestic cropland products available after accounting for the crop-feed requirement of the fixed herd and converted using the total annual human nitrogen requirement. The circle marks the French 2014 reference configuration. The whole explored domain belongs to phase $\mathsf{CE.1.0}$.}
\label{fig:LU-Fsynth}
\end{figure}

Figure~\ref{fig:LU-Fsynth} shows two main effects. First, increasing synthetic fertilization raises $\rho_{POP}$ until cropland production reaches a saturation plateau, after which additional N inputs no longer increase the feeding-capacity indicator. This saturation pattern is consistent with cropland-based feeding-capacity approaches, where yield gains eventually translate into diminishing returns in people nourished per hectare \parencite{cassidy2013redefining,chatzimpiros2023feedingcapacity,harchaoui2019energy}. Second, increasing livestock population lowers $\rho_{POP}$, because a larger herd requires more domestic cropland products for feed, leaving less domestic crop production available for human-equivalent food supply. This livestock effect is visible but less abrupt than the fertilization effect, because additional livestock also increases excreta production and hence recycled N potentially returned to cropland.

At the French 2014 reference point, defined by $f^{synth}_{CROP}=107.2$~kgN\,ha$^{-1}$ and $LU=20.25$~million, the model yields $\rho_{POP}=8.54$~cap\,ha$^{-1}$ cropland, corresponding to a dietary N-equivalent feeding capacity of about 156 million inhabitants. The associated boundary fluxes are $\Phi_{GRP}=2.144\times10^8$~kgN\,yr$^{-1}$ $(23.6$~kgN\,ha$^{-1}$\,yr$^{-1}$),
 $\Phi_{CRP}=-8.446\times10^8$~kgN\,yr$^{-1}$ $(-46.2$~kgN\,ha$^{-1}$\,yr$^{-1}$), and near-zero $\Phi_{LVP}=7.730\times10^6$~kgN\,yr$^{-1}$. With the sign convention used here, the negative $\Phi_{CRP}$ indicates a net export of cropland products at the reference point.

%%%%%%%%%%%%%%%%%%%%%%%%%%
 
\subsubsection{C - Minimal land requirement under dietary change and strict autonomy}
\label{subsection:land-requirement}

In this third case study, we shift from a feeding-capacity perspective to a land-requirement perspective. Rather than asking how many people can be sustained by a given land base, we ask how much agricultural land is minimally required to sustain a fixed population under a prescribed human diet, while imposing strict autonomy in agricultural products. The target population is fixed to the French 2014 reference population, and the two control variables are the mean synthetic nitrogen input to cropland, $f^{synth}_{CROP}$ (kgN\,ha$^{-1}$), and  the vegetal share of the human diet, $\alpha_{POP}$ (N-based).

Unlike in the previous case studies, livestock population is not fixed. It is determined endogenously from the animal share of the human diet, so that the system provides exactly the required amount of livestock products. Cropland and grassland areas are both treated as adjustable variables. For each pair $(f^{synth}_{CROP},\alpha_{POP})$, we determine the minimal cropland and grassland areas compatible with strict product autonomy. The stationary solution is characterized by  zero net boundary fluxes of cropland products, grassland products, and livestock products. The livestock diet is kept at its French reference value, the share of legumes in cropland is fixed at its French reference level, and recycled excreta are redistributed between cropland and grassland proportionally to the corresponding land areas, a rule that holds approximately true for the French reference, see (\ref{eq:EAB}).

\begin{figure}[tbp]
\centering
\includegraphics[width=0.85\textwidth]{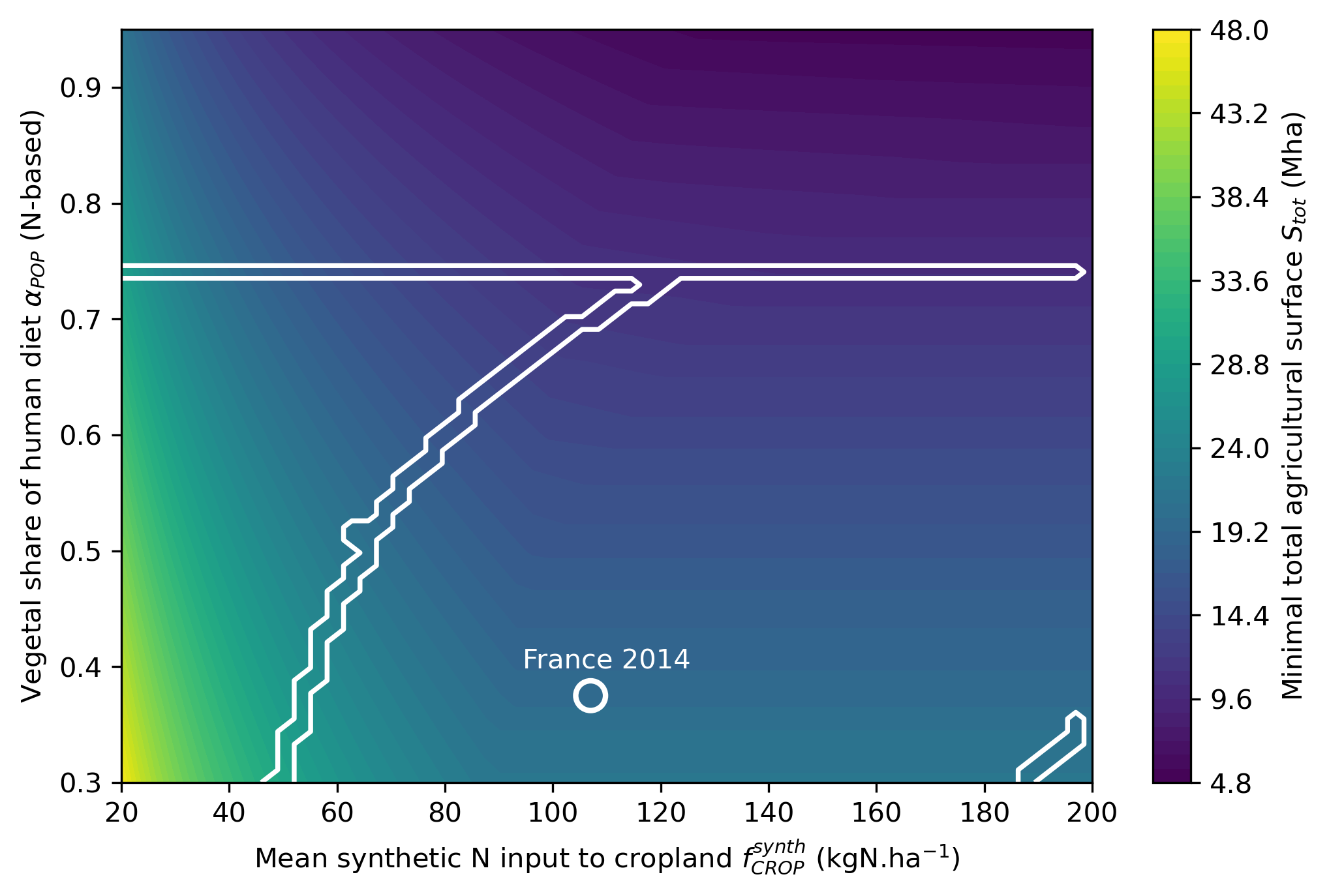}
\caption{Minimal total agricultural surface required to sustain the French reference population under strict product autonomy, as a function of mean synthetic N input to cropland and vegetal share of the human diet. White lines indicate phase boundaries, and the marker denotes the French 2014 configuration. The small lower-right sector corresponds to phase $\mathsf{GE.1}$.}
\label{fig:land-requirement}
\end{figure}

Figure~\ref{fig:land-requirement} shows that the required land area decreases as the human diet becomes more vegetal. Indeed, a more vegetal diet reduces the endogenous livestock population and therefore the land needed to produce animal feed. Also, increasing synthetic fertilization strongly reduces the minimal required area.
%This is consistent with the broader observation that yield gains for major cereals in France and Western Europe have slowed since the 1990s, and that wheat yields in particular have often been described as stagnating over the last decades \parencite{brisson2010wheat,legouis2020wheat}.

At the French 2014 reference point, the model yields a minimal total agricultural surface of about $19.7$~Mha, decomposed into about $10.0$~Mha of cropland and $9.7$~Mha of grassland. This is lower than the approximately $27.4$~Mha of utilised agricultural area reported for France in 2014 \parencite{eurostat2020uaa}. The
difference  should not be interpreted as an agronomic recommendation, but rather as a consequence of the chosen hypotheses (fixed diets,  no export surplus) and coarse-grained representation.

An important result is that the French 2014 diet--fertilization point, under the strict-autonomy closure imposed here, belongs to phase $\mathsf{BE.1}$, instead of $\mathsf{CE.1}$: recycled excreta are distributed evenly between cropland and grassland. The branch $\mathsf{BE.1}$ reflects the livestock-oriented use of cropland products, consistent with the high animal share of the French reference diet.

The phase map helps interpret the contrast between animal-rich, intermediate, and highly vegetal diets. A large part of the diagram lies in $\mathsf{BE.1}$, indicating that under strict autonomy the system is generally organized by a balance between cropland production, grassland production, livestock feeding, and excreta recycling, and characterized by a twofold operational
metabolic cycle, $G_1^{2-LVS}$ (see Fig. 47 in Suppl. Mat.). By contrast, the upper part of the diagram falls into the Trivial regime: when the diet becomes highly vegetarian, the endogenous livestock population becomes small, and internal crop--livestock--excretion loops no longer strongly structure the solution. These configurations may require relatively little land, but they remain dependent on synthetic fertilization or BNF as an external source of reactive nitrogen, with associated environmental and geopolitical concerns \parencite{menegat2022synthetic,alexander2023fertilizer,oecd2024fertiliser}.

At low synthetic inputs, the left part of the diagram enters the crop-excretion subphase $\mathsf{CE.1.0}$ again: this  indicates that the  livestock--excretion--cropland metabolic cycle $G_1^{1-CRP-LVS}$  becomes more structuring when synthetic nitrogen inputs decline. 

 Overall, the figure confirms that reducing synthetic fertilizer use cannot be assessed independently from diet: under strict autonomy, low-input trajectories remain land-efficient only when accompanied by a substantial shift toward more vegetal diets.

\subsection{Broader implications, robustness and limitations}

Recent work has shown that agricultural "feeding capacity" can vary widely depending on diet composition, land use structure, and nitrogen management \parencite{chatzimpiros2023feedingcapacity}. Our contribution is complementary: the reaction-network viewpoint does not only provide indicators, but also a typology of how a territorial agro-food system can operate biophysically. In practice, the analytical phase structure (see Suppl. Mat. for details) provides a translation between (i) familiar agronomic levers (herd size, grazing share, manure routing, fertilizer substitution, wastewater recovery, cropland/grassland endowments) and (ii) the dominant internal recycling loops that actually sustain (or fail to sustain) the system.

\subsubsection{Constrained dynamics}
Most MFA and input-output tools are designed to characterize stationary budgets or to compute one allocation under an objective. A distinctive feature of the present network formalism is that it can also be used to explore time-dependent responses under explicit constraints, without switching to a full numerical model. Concretely, one may keep some quantities fixed because they represent slow or externally imposed constraints (e.g., livestock stock, livestock diet, or a target autonomy constraint), while allowing other degrees of freedom to adjust (e.g., crop allocation between food and feed, boundary fluxes for one product category, or manure routing within an admissible range). This makes it possible to discuss not only feasibility, but also how the system approaches stationarity after a shock, which compartments act as transient bottlenecks, and which constraints can generate non-trivial feedbacks. In a policy context, this provides a mechanism-based way to reason about recovery after disruptions (fertilizer shortage, loss of spreading capacity, changes in grazing share, reduced collection of human excreta) while keeping mass balance and agronomic constraints explicit. 

\noindent Section 10 in Suppl. Mat. presents two cases of stationary states, depending on the imposed constraints: one which is dynamically stable and may illustrate the approach of stationarity after a shock, and one which is dynamically unstable, and causes a total disruption of the nitrogen allocation system.  

\subsubsection{Loss to environment}
Nitrogen losses are not merely an environmental "externality": they are a structural determinant of whether recycling loops can effectively support production. Territorial nitrogen budgets show that losses occur along several segments of the agro-food chain (soils and cropping systems, livestock excretion and manure management, and downstream treatment and dispersion), and that the magnitude and location of these losses strongly condition both agronomic performance and pollution outcomes \parencite{leNoe2017structure,anglade2015soil,billen2012paris}. In our phase language, large losses along the excretion or manure pathway (captured by \textsc{VolRate} and by the effective routing of excreta) tend to "open" the metabolism: recycling becomes too weak to structure the system, which then behaves mainly as a throughput chain sustained by external nitrogen inputs. Conversely, lowering losses (for instance through grazing systems with lower volatilization, improved manure handling, and better synchronization of N returns with crop needs) allows internal loops to become dominant and changes the set of feasible operating modes. This clarifies why "improving nitrogen use efficiency" is not a single scalar objective: the agronomic meaning of reducing losses depends on which loop it enables (crop-livestock coupling, grassland-based ruminant loop, or crop-population recycling loop, the latter being conditional
on a high N recovery from human excreta).

%\subsection{Linking kernel}
%A recurrent conclusion of French socio-metabolic work is that the (re-)coupling between crop and livestock subsystems is a key lever to reduce dependency on industrial inputs and to limit nitrogen leakage, whereas specialization and spatial separation tend to increase feed imports, concentrate surpluses, and displace pollution \parencite{noeBiogeochemicalFunctioningTrajectories,garnier2016reconnecting,anglade2015soil}. In the present framework, this idea is captured by the dominant cyclic "core" (dominant strongly connected component) and by the effective network obtained after collapsing this core. When this kernel is effective, external nitrogen inputs (synthetic fertilizer and biological fixation) are retained longer within the territorial metabolism: nitrogen circulates through crops, livestock and manure return pathways several times before being lost, which reduces the external input required per unit of food produced. When the kernel is weak (poor manure recovery, weak collection, high losses, or allocations that bypass recycling), nitrogen exits the system quickly and production becomes much more dependent on continuous external nitrogen supply.

\subsubsection{Resilience}
In this paper, resilience is used in a practical, biophysical sense: the ability of a territory to absorb perturbations without a qualitative change in its operating mode. When an internal recycling loop is dominant, small variations in management efficiency or routing can often compensate for moderate disturbances (e.g., a temporary increase in food demand can be buffered by a small increase in collection, improved manure routing, or a modest diet shift), because nitrogen can circulate several times within the operational metabolic cycle before being lost. This corresponds to the "multiplier" effect associated with operating cycles described in the Methods and made explicit in the Supplementary Information. When the loop is broken, this buffering disappears: the same small variations can vice versa produce a system-wide reorganization of flows and stocks, or even infeasibility under the imposed constraints. In that sense, phase boundaries behave like biophysical tipping surfaces: close to a boundary, small parameter changes (loss rate, allocation ratio, recovery performance) can change the identity of the dominant loops and therefore of the driving parameters.

\subsubsection{Robustness}
Here robustness refers to a different property: the extent to which the internal state of the operational metabolic cycle is insensitive to external conditions that act outside the core. In agronomic terms, once a dominant recycling kernel is in place, the relative magnitudes of the main core throughputs (e.g., within a crop-livestock-manure loop) are primarily governed by internal allocation structure and conversion efficiencies, and are comparatively less sensitive to external connections that do not belong to the loop (for instance, moderate changes in a peripheral trade flow or a secondary allocation that bypasses the core). This matches the robustness property highlighted in the Methods: core compartments are, to leading order, weakly affected by rates connecting compartments outside the operating cycle. Conversely, in the trivial (no-cycle) regime, there is no protected kernel: the system is directly driven by external inputs and leakages, so changes in fertilizer, trade, or losses propagate without internal buffering. This distinction matters for foresight: two configurations can exhibit similar feeding capacity today, yet differ strongly in robustness because one relies on a dominant internal loop whereas the other is essentially an open chain.

\subsubsection{Parameter sensitivity}
Finally, the phase framework clarifies parameter sensitivity in a way that is directly interpretable for management. In Suppl. Mat, we show that only a small set of levers typically controls whether a given loop becomes dominant: loss rates along manure and excretion management (\textsc{VolRate}), the performance of collection and recovery (\textsc{CollRate} and \textsc{RecRate}), the partition of crop products between food and feed ($k_{CRP\to LVS}$), and the routing of excreta between cropland and grassland. These parameters also differ in pilotability: some are mostly biophysical (yield-response parameters, saturation thresholds), others are constrained yet adjustable via infrastructure and practices (collection, manure handling, grazing share, spreading logistics), and others are primarily societal choices (diet composition, livestock size, crop allocation priorities). As a result, the phase diagrams can be read as a map of which combinations of levers are likely to yield robust circular functioning, which ones produce throughput-driven configurations, and which ones approach feasibility limits.

\subsubsection{Limitations}
This work is intended as a proof of concept and therefore comes with several limitations. First, the model is deliberately aggregated: each compartment represents a broad set of processes (e.g., "cropland products" pools diverse crops and uses), which prevents crop-specific yield responses, feed rations by species, and spatial heterogeneity in soils, climate, and management from being represented explicitly. Second, several coefficients are treated as fixed technical parameters, whereas in reality they vary across regions and farming systems and respond to prices, policies, and learning; our phase diagrams should therefore be read as a qualitative map of biophysical operating modes rather than as quantitative forecasts. Third, yields are represented by a simplified saturating response calibrated on available data; while this captures the existence of a plateau at high nitrogen inputs, it does not represent multi-year soil processes, time-lags in organic nitrogen mineralization, or weather-driven yield variability. Fourth, dynamics are explored only in constrained illustrative settings; a full assessment of transition pathways would require coupling the network formalism to explicit decision rules, investment inertia (e.g., livestock housing, wastewater infrastructure), and spatialized constraints. Finally, we enforce simplified boundary conditions on trade (e.g., exports-only in some case studies) to isolate biophysical mechanisms; real territories are embedded in trade networks, and relaxing these assumptions would change both feasibility domains and the interpretation of "autonomy". These limitations also point to clear extensions: increasing sectoral and spatial resolution, refining yield and loss modules, and integrating behavioral or economic layers while preserving mass-balance consistency and the phase-based interpretability highlighted here.

\section{Conclusion}

We  showed in this work that a phase-based reading of agro-food nitrogen metabolism reveals distinct biophysical operating modes that are not captured by standard input-output or optimization-only approaches. By combining a network representation, analytical phase classification, and three case studies anchored in French reference values, the study identifies how feeding capacity, land requirement, and the balance between recycling and losses depend on the joint configuration of fertilization, crop-livestock coupling, and diet.

The three case studies highlight a consistent result. First, synthetic nitrogen and biological nitrogen fixation are partly substitutable, but only within bounded domains defined by production saturation and the internal structure of nutrient cycling. Second, reducing livestock pressure increases cropland-based feeding capacity, whereas increasing fertilization raises capacity only up to a plateau. Third, strong reductions in synthetic fertilization become difficult to sustain without simultaneous dietary change and redesign of nitrogen sources, because otherwise the land requirement rises sharply or the system becomes strongly dependent on imports. Together, these results show that supply-side and demand-side transitions cannot be assessed independently.

More broadly, the proposed framework provides a transparent way to compare contrasting agrifood configurations through their dominant metabolic phases, rather than only through aggregate performance indicators. This phase-based approach helps connect autonomy, recycling, land use, and vulnerability to external inputs within a single biophysical language. It therefore offers a useful tool for analysing agrifood-system transition pathways and for identifying which configurations are simultaneously productive, autonomous, and structurally robust.

\section*{Declarations}

\subsection*{Funding}  
The authors did not receive support from any organization for the submitted work.

\subsection*{Conflicts of interest}  
The authors have no relevant financial or non-financial interests to disclose.

\subsection*{Acknowledgement}
The authors would like to sincerely thank Julia Le Noë for her insightful discussions regarding this paper.

%\subsection*{Ethics approval}
%Not applicable.

%\subsection*{Consent to participate}  
%Not applicable.

%\subsection*{Consent for publication}  
%Not applicable.

\subsection*{Availability of data and material}  
No standalone dataset was generated for this study. The numerical parameter values used in the model, together with their literature sources and the information needed to interpret them, are reported in this published article and  Suppl. Math.

\subsection*{Code availability}  
Scripts used to reproduce the analyses and generate the figures reported in this article can be provided upon request to the corresponding author.

\subsection*{Authors' contributions}  
AFG: Conceptualization, software, validation, data curation, writing - original draft preparation, visualization.  
JU: Conceptualization, methodology, software, validation, formal analysis, visualization, investigation, writing - original draft preparation.  
LM: Supervision.
OV: Supervision
All authors read and approved the final manuscript.

\subsection*{Use of AI tools}  
During the preparation of this manuscript, the authors used ChatGPT to assist English editing and improve the clarity of the writing. The authors reviewed, edited, and validated all generated text and take full responsibility for the content of the manuscript.

%%%%%%%%%%%%%%%%%%%%%%%%%%%%%%%%%%%%%%%%%%%%%

%\section{Bibliography}
%\bibliographystyle{plainnat}
%\bibliography{bibli}
 
\begingroup
\sloppy
\printbibliography
\endgroup
 
\end{document}